\documentclass[10pt,twoside,reqno]{amsart}
\usepackage{amssymb}
\usepackage{multirow}
\calclayout

\makeatletter

\renewcommand{\@secnumfont}{\bfseries}

\renewcommand{\section}{\@startsection{section}{1}%
  {0mm}{.7\linespacing\@plus\linespacing}{.5\linespacing}
  {\normalfont\bfseries\centering}}

\newcommand{\bibsection}{\@startsection{section}{1}%
  {0mm}{.7\linespacing\@plus\linespacing}{.5\linespacing}
  {\normalfont\scshape\centering}}

\renewcommand{\@biblabel}[1]{#1.}

\newtheorem{thm}{\bf Theorem}[section]

\newtheorem{cor}[thm]{\bf Corollary}

\begin{document}

\vspace{1.3cm}

\title[Some identities of degenerate Euler polynomials]{Some identities of degenerate Euler polynomials associated with degenerate Bernstein polynomials}

\author{Won Joo Kim}
\address{Department of Applied Mathematics, Kyunghee University, Yongin-si, 17104, Republic of Korea}
\email{wjookim@khu.ac.kr}

\author{Dae San Kim}
\address{Department of Mathematics, Sogang University, Seoul 121-742, Republic
of Korea}
\email{dskim@sogang.ac.kr}

\author{Han Young Kim}
\address{Department of Mathematics, Kwnagwoon University, Seoul 139-701, Republic of Korea}
\email{gksdud213@gmail.com}

\author{Taekyun Kim}
\address{Department of Mathematics, Kwangwoon University, Seoul 139-701, Republic
	of Korea}
\email{tkkim@kw.ac.kr}

\subjclass[2010]{11B83, 11S80}
\keywords{degenerate Euler polynomials, degenerate Bernstein polynomials}
\begin{abstract}
In this paper, we investigate some properties and identities for degenerate Euler polynomials in connection with degenerate Bernstein polynomials by means of fermionic $p$-adic integrals on $\Bbb Z_p$ and generating functions. In addition, we study two variable degenerate Bernstein polynomials and the degenerate Bernstein operators.
\end{abstract}
\maketitle

\bigskip
\medskip
\section{Introduction}
Let $p$ be a fixed odd prime number. Throughout this paper, $\Bbb Z_p$, $\Bbb Q_p$ and $\Bbb C_p$, will denote the ring of $p$-adic integers, the field of $p$-adic rational numbers and the completion of algebraic closure of $\Bbb Q_p$, respectively. Let $\nu_p$ be the normalized exponential valuation of $\Bbb C_p$ with $|p \,|_p=p^{-\nu_p (p)}=\frac{1}{p}$. For $\lambda\in\Bbb C_p$ with $|\lambda|_p < p^{-\frac{1}{p-1}}$, the degenerate Euler polynomials are defined by the generating function

\begin{equation}\begin{split}\label{01}
\frac{2}{(1+\lambda t)^{\frac{1}{\lambda}}+1}(1+\lambda t)^{\frac{x}{\lambda}}=\sum_{n=0}^{\infty}{\mathcal E}_{n,\lambda }(x)\frac{t^n}{n!},\ \  (\text{see [1,2]}).
\end{split}\end{equation}
When $x=0$, ${\mathcal E}_{n,\lambda }={\mathcal E}_{n,\lambda }(0)$ are called the degenerate Euler numbers.
The degenerate exponential function is defined by
\begin{equation}\begin{split}\label{02}
e_{\lambda}^{x}(t)=(1+\lambda t)^{\frac{x}{\lambda}}=\sum_{n=0}^{\infty}(x)_{n,\lambda}\frac{t^n}{n!},\ \  (\text{see [6]}),
\end{split}\end{equation}
where
\begin{equation}\begin{split}\label{03}
(x)_{0,\lambda}=1,\  (x)_{n,\lambda}=x(x-\lambda)(x-2\lambda)\cdots(x-(n-1)\lambda), \ \text{for}\ n\geq 1.
\end{split}\end{equation}

From (1), we note that
\begin{equation}\begin{split}\label{04}
{\mathcal E}_{n,\lambda}(x)=\sum_{l=0}^n\binom{n}{l}{\mathcal E}_{l,\lambda}(x)_{n-l,\lambda}, \ (n\geq 0).
\end{split}\end{equation}

Recentely, Kim-Kim introduced the degenerate Bernstein polynomials given by

\begin{equation}\begin{split}\label{05}
\frac{(x)_{k,\lambda}}{k!}t^k (1+\lambda t)^{\frac{1-x}{\lambda}}=\sum_{n=k}^{\infty}B_{k,n}(x|\lambda)\frac{t^n}{n!},\ \  (\text{see [4,5]}).
\end{split}\end{equation}
Thus, by (5), we get

\begin{equation}\begin{split}\label{06}
B_{k,n}(x|\lambda)=\left\{ \begin{array}{lll}\  \binom{n}{k}(x)_{k,\lambda}(1-x)_{n-k,\lambda}, & \mbox{if $n\geq k$},\\\
0, & \mbox{if $n<k$}.  \end{array}  \right. \\
\end{split}\end{equation}
where $n$, $k$ are nonnegative integers.

Let $f$ be a continuous function on $\Bbb Z_p$. Then the degenerate Bernstein operator of order $n$ is given by

\begin{equation}\begin{split}\label{07}
\Bbb B_{n,\lambda}(f|\lambda)&=\sum_{k=0}^{n}f\left(\frac{k}{n}\right)\binom{n}{k}(x)_{k,\lambda}(1-x)_{n-k,\lambda}\\
&=\sum_{k=0}^{n}f\left(\frac{k}{n}\right)B_{k,n}(x|\lambda),\ \  (\text{see [4,5]}).
\end{split}\end{equation}
The fermionic $p$-adic integral on $\Bbb Z_p$ is defined by Kim as

\begin{equation}\begin{split}\label{08}
\int_{\Bbb Z_p } f(x) d\mu_{-1} (x)=  \lim_{N \rightarrow \infty} \sum_{x=0}^{p^N-1} f(x)(-1)^x ,\ \  (\text{see [3,5]}).
\end{split}\end{equation}
By (8), we get
\begin{equation}\begin{split}\label{09}
\int_{\Bbb Z_p } f(x+1) d\mu_{-1} (x)+\int_{\Bbb Z_p } f(x) d\mu_{-1} (x)=2f(0),\ \  (\text{see [3,7-12]}).
\end{split}\end{equation}
From (8), we note that
\begin{equation}\begin{split}\label{10}
\int_{\Bbb Z_p }(1+\lambda t)^{\frac{x+y}{\lambda}}d\mu_{-1}(y)=\frac{2}{(1+\lambda t)^{\frac{1}{\lambda}}+1}(1+\lambda t)^{\frac{x}{\lambda}}=\sum_{n=0}^{\infty}{\mathcal E}_{n,\lambda}(x)\frac{t^n}{n!}.
\end{split}\end{equation}
On the other hand,
\begin{equation}\begin{split}\label{11}
\int_{\Bbb Z_p }(1+\lambda t)^{\frac{x+y}{\lambda}}d\mu_{-1}(y)=\sum_{n=0}^\infty\int_{\Bbb Z_p }(x+y)_{n,\lambda}d\mu_{-1}(y)\frac{t^n}{n!}.
\end{split}\end{equation}
By (10) and (11), we get
\begin{equation}\begin{split}\label{12}
\int_{\Bbb Z_p }(x+y)_{n,\lambda}d\mu_{-1}(y)={\mathcal E}_{n,\lambda}(x),\ \ (n\geq 0),\ \  (\text{see [4,5]}).
\end{split}\end{equation}
In this paper, we investigate some properties and identities for degenerate Euler polynomials in connection with degenerate Bernstein polynomials by means of fermionic $p$-adic integrals on $\Bbb Z_p$ and generating functions. In addition, we study two variable degenerate Bernstein polynomials and the degenerate Bernstein operators.

\section{Degenerate Euler and Bernstein polynomials}
From (1), we note that
\begin{equation}\begin{split}\label{13}
2&=\sum_{n=0}^\infty\left(\sum_{m=0}^{n}\binom{n}{m}{\mathcal E}_{m,\lambda}(1)_{n-m,\lambda}+{\mathcal E}_{n,\lambda}\right)\frac{t^n}{n!}\\
 &=\sum_{n=0}^\infty\left({\mathcal E}_{n,\lambda}(1)+{\mathcal E}_{n,\lambda}\right)\frac{t^n}{n!}.
\end{split}\end{equation}
Comparing the coefficients on both sides of (13), we have
\begin{equation}\begin{split}\label{14}
{\mathcal E}_{n,\lambda}(1)+{\mathcal E}_{n,\lambda}=2\delta_{0,n},\  (n,k\geq 0),
\end{split}\end{equation}
where $\delta_{n,k}$ is the Kronecker's symbol.

By (1), we easily get
\begin{equation}\begin{split}\label{15}
{\mathcal E}_{n,\lambda}(1-x)=(-1)^n {\mathcal E}_{n,-\lambda}(x),\ \ (n\geq 0).
\end{split}\end{equation}
From (1), (4) and (14), we note that
\begin{equation}\begin{split}\label{16}
{\mathcal E}_{n,\lambda}(2)&=\sum_{l=0}^{n}\binom{n}{l}{\mathcal E}_{l,\lambda}(1)(1)_{n-l,\lambda}\\
&=(1)_{n,\lambda}+\sum_{l=1}^{n}\binom{n}{l}{\mathcal E}_{l,\lambda}(1)(1)_{n-l,\lambda}\\
&=2(1)_{n,\lambda}-\sum_{l=0}^{n}\binom{n}{l}(1)_{n-l,\lambda}{\mathcal E}_{l,\lambda}\\
&=2(1)_{n,\lambda}+{\mathcal E}_{n,\lambda},
\end{split}\end{equation}
where $n$ is a positive integer.

Therefore, by (16), we obtain the following theorem.
\begin{thm}
For $n\in\Bbb N$, we have
\begin{equation*}\begin{split}\label{2.1}
 {\mathcal E}_{n,\lambda}(2)=2(1)_{n,\lambda}+{\mathcal E}_{n,\lambda}.
\end{split}\end{equation*}
\end{thm}
Note that
\begin{equation}\begin{split}\label{17}
(1-x)_{n,\lambda}=(-1)^n (x-1)_{n,-\lambda},\,\, (n \geq 0).
\end{split}\end{equation}
Therefore, by (12), (15) and (17), we easily get
\begin{equation}\begin{split}\label{18}
\int_{\Bbb Z_p }(1-x)_{n,\lambda}d\mu_{-1}(x)&=(-1)^n \int_{\Bbb Z_p }(x-1)_{n,-\lambda}d\mu_{-1}(x)\\
&=\int_{\Bbb Z_p }(x+2)_{n,\lambda}d\mu_{-1}(x).
\end{split}\end{equation}
Therefore, by (18) and Theorem 2.1, we obtain the following theorem.
\begin{thm}
For $n\in\Bbb N$, we have
\begin{equation*}\begin{split}\label{2}
 \int_{\Bbb Z_p }(1-x)_{n,\lambda}d\mu_{-1}(x)=\int_{\Bbb Z_p }(x+2)_{n,\lambda}d\mu_{-1}(x)=2(1)_{n,\lambda}+\int_{\Bbb Z_p }(x)_{n,\lambda}d\mu_{-1}(x).
\end{split}\end{equation*}
\end{thm}

\begin{cor}
For $n\in\Bbb N$, we have
\begin{equation*}\begin{split}\label{2.3}
(-1)^{n} {\mathcal E}_{n,-\lambda}(-1)=2(1)_{n,\lambda}+{\mathcal E}_{n,\lambda}={\mathcal E}_{n,\lambda}(2).
\end{split}\end{equation*}
\end{cor}

By (4), we get
\begin{equation}\begin{split}\label{19}
{\mathcal E}_{n,\lambda}(1-x)&=\sum_{l=0}^{n}\binom{n}{l}(1-x)_{n-l,\lambda}{\mathcal E}_{l,\lambda}\\
&=\sum_{l=0}^{n}\binom{n}{l}(x)_{l,\lambda}(1-x)_{n-l,\lambda}\frac{{\mathcal E}_{l,\lambda}}{(x)_{l,\lambda}}\\
&=\sum_{l=0}^{n}B_{l,n}(x|\lambda){\mathcal E}_{l,\lambda}\frac{1}{(x)_{l,\lambda}}.
\end{split}\end{equation}

Let
\begin{equation}\begin{split}\label{20}
\frac{1}{(x)_{l,\lambda}}=\frac{1}{x(x-\lambda)(x-2\lambda)\cdots (x-(l-1)\lambda)}=\sum_{k=0}^{l-1}\frac{A_k}{x-k\lambda},\ (l\in\Bbb N).
\end{split}\end{equation}

Then we have
\begin{equation}\begin{split}\label{21}
A_k=\lambda^{1-l}\prod^{l-1}_{\substack{i=0,\\i\neq k}} \left(\frac{1}{k-i}\right)=\lambda^{1-l}\frac{(-1)^{k-l-1}}{k!(l-1-k)!}=\frac{\lambda^{1-l}}{(l-1)!}\binom{l-1}{k}(-1)^{k-l-1}.
\end{split}\end{equation}
By (20) and (21), we get
\begin{equation}\begin{split}\label{22}
A_k=\frac{(-\lambda)^{1-l}}{(l-1)!}\binom{l-1}{k}(-1)^{k}.
\end{split}\end{equation}
From (20) and (22), we have
\begin{equation}\begin{split}\label{23}
\frac{1}{(x)_{l,\lambda}}=\sum_{k=0}^{l-1}\frac{(-1)^k}{(l-1)!}\binom{l-1}{k}\frac{(-\lambda)^{1-l}}{x-k\lambda},\ (l\in\Bbb N).
\end{split}\end{equation}
By (19) and (20), we get
\begin{equation}\begin{split}\label{24}
{\mathcal E}_{n,\lambda}(1-x)&=\sum_{l=0}^{n}B_{l,n}(x|\lambda){\mathcal E}_{l,\lambda}\frac{1}{(x)_{l,\lambda}}\\
&=(1-x)_{n,\lambda}+\sum_{l=1}^{n}B_{l,n}(x|\lambda){\mathcal E}_{l,\lambda}\frac{1}{(x)_{l,\lambda}}\\
&=(1-x)_{n,\lambda}+\sum_{l=1}^{n}B_{l,n}(x|\lambda){\mathcal E}_{l,\lambda}\frac{(-\lambda)^{1-l}}{(l-1)!}\sum_{k=0}^{l-1}(-1)^k\binom{l-1}{k}\frac{1}{x-k\lambda}.
\end{split}\end{equation}
Therefore, by (24), we obtain the following theorem.

\begin{thm}
For $n\geq 0$, we have
\begin{equation*}\begin{split}\label{2.1}
 {\mathcal E}_{n,\lambda}(1-x)=(1-x)_{n,\lambda}+\sum_{l=1}^{n}B_{l,n}(x|\lambda){\mathcal E}_{l,\lambda}\frac{(-\lambda)^{1-l}}{(l-1)!}\sum_{k=0}^{l-1}(-1)^k\binom{l-1}{k}\frac{1}{x-k\lambda}.
\end{split}\end{equation*}
\end{thm}

\begin{cor}
For $n\geq 0$, we have
\begin{equation*}\begin{split}\label{2.1}
{\mathcal E}_{n,\lambda}(2)=(2)_{n,\lambda}- \sum_{l=1}^{n}B_{l,n}(-1|\lambda){\mathcal E}_{l,\lambda}\frac{(-\lambda)^{1-l}}{(l-1)!}\sum_{k=0}^{l-1}(-1)^{k}\binom{l-1}{k}\frac{1}{1+k\lambda}.
\end{split}\end{equation*}
\end{cor}

For $k\in\Bbb N$, the higher order degenerate Euler polynomials are given by the generating function
\begin{equation}\begin{split}\label{25}
\left(\frac{2}{(1+\lambda t)^{\frac{1}{\lambda}}+1}\right)^k(1+\lambda t)^{\frac{x}{\lambda}}=\sum_{n=0}^{\infty}{\mathcal E}_{n,\lambda}^{(k)}(x)\frac{t^n}{n!}.
\end{split}\end{equation}
From (5) and (25), we note that
\begin{equation}\begin{split}\label{26}
\sum_{n=0}^{\infty}\frac{1}{\binom{n+k}{n}}B_{k,n+k}(x|\lambda)\frac{t^n}{n!}&=(x)_{k,\lambda}(1+\lambda t)^{\frac{1-x}{\lambda}}\\
&=\frac{(x)_{k,\lambda}}{2^k}\sum_{l=0}^{k}\binom{k}{l}\left(\frac{2}{(1+\lambda t)^{\frac{1}{\lambda}}+1}\right)^k(1+\lambda t)^{\frac{1-x+l}{\lambda}}\\
&=\frac{(x)_{k,\lambda}}{2^k}\sum_{n=0}^{\infty}\left(\sum_{l=0}^{k}\binom{k}{l}{\mathcal E}_{n,\lambda}^{(k)}(1-x+l)\right)\frac{t^n}{n!}.
\end{split}\end{equation}
Therefore, by comparing the cofficients on both sides of (26), we obtain the following theorem.

\begin{thm}
For $n,k\in\Bbb N$, we have
\begin{equation*}\begin{split}\label{2.1}
\frac{2^k}{\binom{n+k}{n}}B_{k,n+k}(x|\lambda)=(x)_{k,\lambda}\sum_{l=0}^{k}\binom{k}{l}{\mathcal E}_{n,\lambda}^{(k)}(1-x+l).
\end{split}\end{equation*}
\end{thm}
Let $f$ be a continuous function on $\Bbb Z_p$. For $x_1,x_2\in\Bbb Z_p$, we consider the degenerate Bernstein operator of order $n$ given by
\begin{equation}\begin{split}\label{27}
\Bbb B_{n,\lambda}(f|x_{1},x_{2})&=\sum_{k=0}^{n}f\left(\frac{k}{n}\right)\binom{n}{k}(x_1)_{k,\lambda}(1-x_2)_{n-k,\lambda}\\
&=\sum_{k=0}^{n}f\left(\frac{k}{n}\right)B_{k,n}(x_{1},x_{2}|\lambda),
\end{split}\end{equation}
where
\begin{equation}\begin{split}\label{28}
B_{k,n}(x_{1},x_{2}|\lambda)=\binom{n}{k}(x_{1})_{k,\lambda}(1-x_2)_{n-k,\lambda},
\end{split}\end{equation}
where $n,k$ are nonnegative integers.

Here, $B_{k,n}(x_{1},x_{2}|\lambda)$ are called two variable degenerate Bernstein polynomials of degree $n$.

From (28), we note that
\begin{equation}\begin{split}\label{29}
\sum_{n=k}^{\infty}B_{k,n}(x_{1},x_{2}|\lambda)\frac{t^n}{n!}&=\sum_{n=k}^{\infty}\binom{n}{k}(x_{1})_{k,\lambda}(1-x_2)_{n-k,\lambda}\frac{t^n}{n!}\\
&=\sum_{n=k}^{\infty}\frac{(x_{1})_{k,\lambda}(1-x_2)_{n-k,\lambda}}{k!(n-k)!}t^n\\
&=\frac{(x_{1})_{k,\lambda}}{k!}t^k\sum_{n=0}^{\infty}\frac{(1-x_2)_{n,\lambda}}{n!}t^n\\
&=\frac{(x_{1})_{k,\lambda}}{k!}t^k(1+\lambda t)^{\frac{1-x_2}{\lambda}}\\
&=\frac{(x_{1})_{k,\lambda}}{k!}t^ke_{\lambda}^{1-x_2}(t).
\end{split}\end{equation}
Thus, by (29), we get
\begin{equation}\begin{split}\label{30}
\frac{(x_{1})_{k,\lambda}}{k!}t^k(1+\lambda t)^{\frac{1-x_2}{\lambda}}=\sum_{n=k}^{\infty}B_{k,n}(x_1,x_2 |\lambda)\frac{t^n}{n!},
\end{split}\end{equation}
where $k$ is a nonnegative integer.
By (28), we easily get
\begin{equation}\begin{split}\label{32}
B_{k,n}(x_1,x_2|\lambda)&=\binom{n}{k}(1-(1-x_1))_{n-(n-k),\lambda}(1-x_2)_{n-k,\lambda}\\
&=B_{n-k,n}(1-x_2,1-x_1 |\lambda).
\end{split}\end{equation}
Now, we observe that
\begin{equation}\begin{split}\label{33}
&(1-x_2-(n-k-1)\lambda)B_{k,n-1}(x_1,x_2|\lambda)+(x_1-(k-1)\lambda)B_{k-1,n-1
}(x_1,x_2 |\lambda)\\
&=(1-x_2-(n-k-1)\lambda)\binom{n-1}{k}(x_1)_{k,\lambda}(1-x_2)_{n-1-k,\lambda}\\&+(x_1-(k-1)\lambda)\binom{n-1}{k-1}(x_1)_{k-1,\lambda}(1-x_2)_{n-k,\lambda}\\
&=\binom{n}{k}(x_1)_{k,\lambda}(1-x_2)_{n-k,\lambda}
=B_{k,n}(x_1,x_2|\lambda),\ \ (n,k\in\Bbb N).
\end{split}\end{equation}
Therefore, by (32), we obtain the following theorem.
\begin{thm}
For $n,k\in\Bbb N$, we have
\begin{equation*}\begin{split}\label{2.1}
&(1-x_2-(n-k-1)\lambda)B_{k,n-1}(x_1,x_2|\lambda)+(x_1-(k-1)\lambda)B_{k-1,n-1
}(x_1,x_2 |\lambda)\\&=B_{k,n}(x_1,x_2|\lambda).
\end{split}\end{equation*}
\end{thm}
If $f=1$, then we have from (27),
\begin{equation}\begin{split}\label{34}
\Bbb B_{n,\lambda}(1|x_1,x_2)&=\sum_{k=0}^{n}B_{k,n}(x_1,x_2|\lambda)
=\sum_{k=0}^{n}\binom{n}{k}(x_1)_{k,\lambda}(1-x_2)_{n-k,\lambda}\\
&=(1+x_1-x_2)_{n,\lambda}.
\end{split}\end{equation}
If $f(t)=t$, then we also get from (27) that for $n\in\Bbb N$ and $x_1,x_2\in\Bbb Z_p$,
\begin{equation}\begin{split}\label{35}
\Bbb B_{n,\lambda}(t|x_1,x_2)&=\sum_{k=0}^{n}\frac{k}{n}\binom{n}{k}(x_1)_{k,\lambda}(1-x_2)_{n-k,\lambda}\\
&=(x_1)_{1,\lambda}(x_1+1-\lambda-x_2)_{n-1,\lambda}.
\end{split}\end{equation}
Hence,
\begin{equation}\begin{split}\label{36}
(x_1)_{1,\lambda}=\frac{1}{(x_1+1-\lambda-x_2)_{n-1,\lambda}}\Bbb B_{n,\lambda}(t|x_1,x_2).
\end{split}\end{equation}
By the same method, we get
\begin{equation}\begin{split}\label{37}
&\Bbb B_{n,\lambda}(t^2|x_1,x_2)\\
&=\frac{1}{n}(x_1)_{1,\lambda}(1+x_1-\lambda-x_2)_{n-1,\lambda}+\frac{n-1}{n}(x_1)_{2,\lambda}(1+x_1-2\lambda-x_2)_{n-2,\lambda}.
\end{split}\end{equation}
Note that
\begin{equation*}\begin{split}
\lim_{n\rightarrow \infty}\Big(\lim_{\lambda\rightarrow 0}\Bbb B_{n,\lambda}(t^2|x,x)\Big)  = \lim_{n\rightarrow \infty }\left(\frac{x}{n}+\frac{n-1}{n}x^2\right)=x^2.
\end{split}\end{equation*}
Now, we observe that
\begin{equation}\begin{split}\label{38}
\sum_{k=1}^{n}\frac{\binom{k}{1}}{\binom{n}{1}}B_{k,n}(x_1,x_2|\lambda)&=\sum_{k=1}^{n}\binom{n-1}{k-1}(x_1)_{k,\lambda}(1-x_2)_{n-k,\lambda}\\
&=\sum_{k=0}^{n-1}\binom{n-1}{k}(x_1)_{k+1,\lambda}(1-x_2)_{n-1-k,\lambda}\\
&=(x_1)_{1,\lambda}(x_1+1-\lambda-x_2)_{n-1,\lambda}.
\end{split}\end{equation}
Thus, by (37), we get
\begin{equation*}\begin{split}
(x_1)_{1,\lambda}=\frac{1}{(1+x_1-x_2-\lambda)}_{n-1,\lambda}\sum_{k=1}^{n}\frac{\binom{k}{1}}{\binom{n}{1}}B_{k,n}(x_1,x_2|\lambda).
\end{split}\end{equation*}
By the same method, we get
\begin{equation*}\begin{split}
(x_1)_{2,\lambda}=\frac{1}{(1+x_1-x_2-2\lambda)}_{n-2,\lambda}\sum_{k=2}^{n}\frac{\binom{k}{2}}{\binom{n}{2}}B_{k,n}(x_1,x_2|\lambda).
\end{split}\end{equation*}
Continuing this process, we have
\begin{equation}\begin{split}\label{39}
(x_1)_{i,\lambda}=\frac{1}{(1+x_1-x_2-i\lambda)}_{n-i,\lambda}\sum_{k=i}^{n}\frac{\binom{k}{i}}{\binom{n}{i}}B_{k,n}(x_1,x_2|\lambda),\ \ (i\in\Bbb N).
\end{split}\end{equation}
\begin{thm}
For $i\in\Bbb N$, we have
\begin{equation*}\begin{split}\label{7}
(x)_{i,\lambda}=\frac{1}{(1+x_1-x_2-i\lambda)}_{n-i,\lambda}\sum_{k=i}^{n}\frac{\binom{k}{i}}{\binom{n}{i}}B_{k,n}(x_1,x_2|\lambda),\ \ (i\in\Bbb N).
\end{split}\end{equation*}
\end{thm}
Taking double fermionic $p$-adic integral on $\Bbb Z_p$, we get the following equation:
\begin{equation}\begin{split}\label{41}
&\int_{\Bbb Z_p } \int_{\Bbb Z_p } B_{k,n}(x_1,x_2|\lambda)d\mu_{-1}(x_1)d\mu_{-1}(x_2)\\&=\binom{n}{k}\int_{\Bbb Z_p }(x_1)_{k,\lambda}d\mu_{-1}(x_1)\int_{\Bbb Z_p }(1-x_2)_{n-k,\lambda}d\mu_{-1}(x_2).
\end{split}\end{equation}
Therefore, by (39) and Theorem 2.2, we obtain the following theorem.
\begin{thm}\label{thm9}
For $n,k\geq 0$, we have
\begin{equation*}\begin{split}\label{8}
&\int_{\Bbb Z_p } \int_{\Bbb Z_p } B_{k,n}(x_1,x_2|\lambda)d\mu_{-1}(x_1)d\mu_{-1}(x_2)\\
&\quad\quad\quad=\left\{ \begin{array}{lll}\  \binom{n}{k}\mathcal{E}_{k,\lambda}\left(2(1)_{n-k,\lambda}+\mathcal{E}_{n-k,\lambda}\right), & \mbox{if $n> k$},\\\
\mathcal{E}_{n,\lambda}, & \mbox{if $n=k$}.  \end{array}  \right. \\
\end{split}\end{equation*}
\end{thm}
We get from the symmetric properties of two variable degenerate Bernstein polynomials that for $n,k\in \Bbb N$ with $n>k$,
\begin{equation}\begin{split}\label{42}
&\int_{\Bbb Z_p } \int_{\Bbb Z_p }B_{k,n}(x_1,x_2|\lambda)d\mu_{-1}(x_1)d\mu_{-1}(x_2)
\\&=\sum_{l=0}^{k}\binom{n}{k}\binom{k}{l}(-1)^{k+l}(1)_{l,\lambda}\\
&\quad\quad\times \int_{\Bbb Z_p}
\int_{\Bbb Z_p }(1-x_1)_{k-l,-\lambda}(1-x_2)_{n-k,\lambda}d\mu_{-1}(x_1)d\mu_{-1}(x_2)\\
&=\binom{n}{k}\int_{\Bbb Z_p }(1-x_2)_{n-k,\lambda}d\mu_{-1}(x_2)\left\{ (1)_{k,\lambda}+\sum_{l=0}^{k-1}\binom{k}{l}(-1)^{k+l}(1)_{l,\lambda} \mathcal{E}_{k-l,-\lambda} (2)\right\} \\
&=\binom{n}{k}\mathcal{E}_{n-k,\lambda}(2)\left\{(1)_{k,\lambda}+\sum_{l=0}^{k-1}\binom{k}{l}(-1)^{k+l}(1)_{l,\lambda}\mathcal{E}_{k-l,-\lambda}(2)\right\}.
\end{split}\end{equation}
Therefore, by Theorem \ref{thm9} and (40), we obtain the following theorem.
\begin{thm}
For $k\in\Bbb N$, we have
\begin{equation*}\begin{split}\label{10}
\mathcal{E}_{k,\lambda}=(1)_{k,\lambda}+\sum_{l=0}^{k-1}\binom{k}{l}(-1)^{k+l}(1)_{l,\lambda}\left(\mathcal{E}_{k-l,-\lambda}+2(1)_{k-l,-\lambda}\right).
\end{split}\end{equation*}
\end{thm}
Note that
\begin{equation*}\begin{split}\label{39}
&\sum_{l=0}^{k-1}\binom{k}{l}(-1)^{k+l}(1)_{l,\lambda}(1)_{k-l,-\lambda}\\
&=(-1)^k\Big(\sum_{l=0}^{k}\binom{k}{l}(-1)_{l,-\lambda}(1)_{k-l,-\lambda}-(-1)_{k,-\lambda}\Big)\\
&=(-1)^k\Big((0)_{k,-\lambda}-(-1)_{k,-\lambda}\Big)\\
&=-(1)_{k,\lambda}.
\end{split}\end{equation*}

\begin{cor}
For $k\in\Bbb N$, we have
\begin{equation*}\begin{split}\label{10}
\mathcal{E}_{k,\lambda}=-(1)_{k,\lambda}+\sum_{l=0}^{k-1}\binom{k}{l}(-1)^{k+l}(1)_{l,\lambda}\mathcal{E}_{k-l,-\lambda}.
\end{split}\end{equation*}
\end{cor}

\section{Conclusions}

\indent In [1,2], Carlitz initiated the study of degenerate versions of some special polynomials and numbers, namely the degenerate Bernoulli and Euler polynomials and numbers. Here we would like to draw the attention of the readers that T. Kim and his colleagues have been introducing various degenerate polynomials and numbers and investigating their properties, some identities related to them and their applications by means of generating functions, combinatorial methods, umbral calculus, $p$-adic analysis and differential equations (see [4,5] and the references therein). It is amusing that this line of study led them even to the introduction of degenerate gamma functions and degenerate Laplace transforms (see [7]). These already demonstrate that studying various degenerate versions of known special numbers and polynomials can be very promising and rewarding. Furthermore, we can hope that many applications will be found not only in mathematics but also in sciences and engineering. \\
\indent In this paper, we investigated some properties and identities for degenerate Euler polynomials in connection with degenerate Bernstein polynomials and operators which were recently introduced as degenerate versions of the classical Bernstein polynomials and operators. This has been done by means of fermionic $p$-adic integrals on $\Bbb Z_p$ and generating functions. In addition, we studied two variable degenerate Bernstein polynomials and the degenerate Bernstein operators.

\vspace{0.1in}

{\bf Acknowledgements}

 The fourth author’s work in this paper was conducted during the sabbatical year of Kwangwoon
University in 2018.

\vspace{0.1in}

{\bf Competing interests}

The authors declare that they have no competing interests.

\vspace{0.1in}

{\bf Authors’ contributions}

All authors contributed equally to the manuscript and typed, read, and approved the final manuscript.

\vspace{0.1in}

{\bf Funding}

This research received no external funding.

\end{document}